\def\IH{{\Bbb H}} 
 
\def\IR{{\Bbb R}} 
 
\def\IS{{\Bbb S}} 
 
\def\IK{{\Bbb K}}

\def\IC{\Bbb C} 
\def\oC{\hat{\IC}} 
\def\ID{{\Bbb D}}
\def\oD{\overline{{\Bbb D}}}

\def\zbar{{\overline{z}}} 
 
\def\wbar{{\overline{w}}}

\def\mA{{\cal A}}
\def\mE{{\mathsf{E}}}
\def\mF{{\mathsf{F}}}
\def\mH{{\mathsf{H}}}

\documentclass[12pt]{article} 

\usepackage[psamsfonts]{amssymb} 
\usepackage{amsfonts,amsmath} 
\usepackage{enumitem}

\usepackage{epsfig,multicol}

\DeclareMathOperator{\varg}{Varg}
\DeclareMathOperator{\diam}{diam}

\newtheorem{theorem}{Theorem}
\newtheorem{lemma}{Lemma}

\newtheorem{conjecture}{Conjecture}[section] 
 
\numberwithin{equation}{section}

\title{Diffeomorphic solutions of Ahlfors-Hopf equations}
\author{Gaven Martin \& Cong Yao\thanks{\tiny
Work of both authors partially supported by the New Zealand Marsden Fund. C. Yao is partially supported by the Young Scientist Program of the Ministry of Science and Technology of China (No. 2021YFA1002200), the National Natural Science Foundation of China (No. 12401096, No. 12101362), the Natural Science Foundation of Shandong Province (No. ZR2024QA035, No. ZR2022YQ01), the Fundamental Research Funds for the Central Universities.
\newline
Institute for Advanced Study, 
Massey University,  Auckland,
New Zealand.
\newline
email: g.j.martin@massey.ac.nz  \newline
Research Center for Mathematics and Interdisciplinary Sciences, Shandong University, 266237, Qingdao and Frontiers Science Center for Nonlinear Expectations, Ministry of Education, P. R. China. \newline
email:  c.yao@sdu.edu.cn 
\newline
\newline
{\bf MSC Subject.}  30C62 31A05 49J10  } }
\date{}
\begin{document} 
\maketitle

\begin{abstract}
{\footnotesize
We consider the boundary the value problem for extremal functions of $L^p$-mean distortion and the associated Teichm\"uller spaces, interpolating between classical extremal quasiconformal mapping,  and recent approachs through harmonic mappings (of extreme Dirichlet energy). This directly leads to the $L^p$-Teichm\"uller theory of Riemann surfaces of finite area.  
Here, given convex increasing $\mA$, we focus on the Alhfors-Hopf differential
\vskip-9pt
\[  
\Phi=\mA(\IK(w,h))h_w\,\overline{h_\wbar}\, \eta(h),
\]
\vskip-4pt
where $h=f^{-1}$ is the pseudo-inverse of an extremal mapping  $f$ for 
\vskip-9pt
\[ \inf\;\int_\ID \mA(\IK(z,f)) \; dz, \quad\quad \IK(z,f) =   \frac{|f_z|^2+|f_\zbar|^2}{|f_z|^2-|f_\zbar|^2}. \]
where the infimum is taken over those homeomorphisms of finite distortion $f:\overline{\ID}\to\overline{\ID}$ with $f|\IS=f_0$,  typically a quasisymmetric barrier.  The inner-variational equations, an analogue of the Euler-Lagrange equations, show $\Phi$ is holomorphic at an extremal.  Exploiting this Ahlfors-Hopf differential, we prove that an extreme point $f$ is a local diffeomorphism in $\ID$, resolving some conjectures in \cite{MY3}. Then, given a pair of closed hyperbolic Riemann surfaces $(\Sigma_i,\sigma_i)$ and $1\leq p<\infty$,  we show there is a unique diffeomorphic minimiser to the problem
\vskip-12pt
\[ \min_{f\in \, {\rm Diff}(\Sigma_1,\Sigma_2)} \int_{\Sigma_1} \IK^p(w,f) \;  d\sigma_1(w)\]
As $p\to\infty$ these maps converge to the quasiconformal extremal, but it is not a diffeomorphism,  smoothness is lost on the divisor of an associated holomorphic quadratic differential.}
\end{abstract}

\section{Introduction} 
In \cite{AIMO} and \cite{IKO} the authors posed the following conjecture concerning the boundary value problem for self mappings of the unit disk $\ID$ with finite distortion.  Suppose $\mA:[1,\infty)\to [1,\infty)$ is convex increasing.  Define the $\mA$-conformal energy of a mapping $f:\ID\to\ID$ of finite distortion by
\[ \mE_{\mA}(f):=\int_\ID\mA(\IK(z,f))\; dz, \quad \IK(z,f) = \frac{|f_z|^2+|f_\zbar|^2}{|f_z|^2-|f_\zbar|^2}. \]
Here a mapping is of finite distortion if $f\in W^{1,1}_{loc}(\ID)$, the Jacobian determinant $J(z,f)=|f_z|^2-|f_\zbar|^2\in L^{1}_{loc}(\ID)$ and $\IK(z,f)$ is finite almost everywhere.  The basics of the theory of mappings of finite distortion can be found in the monograph \cite{AIM}.   Basically,  they are linked to degenerate elliptic Beltrami equations in the same way that quasiconformal mappings are linked to uniformly elliptic Beltrami equations. 
\begin{conjecture}
Suppose $f_0:\ID\to\ID$ is a prescribed barrier function with $\mE_{\mA}(f_0)<\infty$.  The variational problem  
\begin{equation}\label{1.1}
\inf_{f} \big\{ \mE_{\mA}(f),\quad \mbox{$f$ has finite distortion and $f|_\IS=f_0|_\IS$} \big\}
\end{equation}
 has a unique diffeomorphic minimiser.
 \end{conjecture}
 The most interesting cases are $\mA(t)=t^p$, $p\geq 1$,  or $\mA(t)=e^{pt}$, $p>0$.  We generally assume that the barrier $f_0$ is homeomorphic, or even quasiconformal.  In the case of $L^1$-means of distortion ($\mA(t)=t$)  this conjecture is true,  the harmonic extension of $(f_0|_\IS)^{-1}$ is the unique diffeomorphic minimiser, \cite{AIMO,IO}.  Further,  as $p\to \infty$ any sequence of minimisers $f_p$,  whether smooth or not, converges to an extremal quasiconformal mapping (this is basically due to Ahlfors \cite{Ah3} in his celebrated proof of Teichm\"uller's theorem).  Note that this extremal quasiconformal mapping will most often {\em not} be smooth.  Other interesting extremal problems are studied in \cite{3},  ultimately leading to the solution of the Nitsche conjecture \cite{6}, the Gr\"otzsch problem  \cite{9} and Teichm\"uller's problem, \cite{2,TP}. 
 
\medskip

In this article we consider a slightly more  general energy functional by including a weight $\eta\in C^\infty(\ID)$, $\eta\geq 1$.  The problem is to minimise the functional
\begin{equation}\label{1.2}
\mE_{\mA,\eta}(f):=\int_\ID\mA(\IK(z,f))\; \eta(z)\, dz,  
\end{equation}
where $f:\overline{\ID}\to\overline{\ID}$ is a homeomorphism of finite distortion, with prescribed boundary values $f|_\IS=f_0|_\IS$. We may relax the assumption $f$ is a homeomorphism to $f$ is a surjection. Note that as $\mA\geq 1$ we must have $\eta\in L^1(\ID)$.

Inner variation,  see \cite{AIMO}, leads to the following distributional Euler-Lagrange equations on the ``$f$-side'':
\begin{equation}\label{1.3}
2\int_\ID\mA'(\IK(z,f))\frac{\overline{\mu_f}}{1-|\mu_f|^2}\eta\varphi_\zbar=\int_\ID\mA(\IK(z,f))(\eta\varphi)_z,\quad\forall\varphi\in C_0^\infty(\ID).
\end{equation}
where 
\[ \mu_f(z) = \frac{f_\zbar(z)}{f_z(z)}, \quad \mbox{the Beltrami coefficient of $f$}.\]
This equation holds for any extremal $f$.  However if we consider the {\em inverse functional}
\begin{equation}\label{1.4}
\mathsf{E}_{\mA,\eta}^\ast(h):=\int_\ID\mA(\IK(w,h))\;J(w,h)\; \eta(h)\; dw,  
\end{equation}
where $h$ has finite distortion and  $h|_\IS=(f_0)^{-1}|_\IS$, and obtained by the change of variables $z=h(w)$ ($f(z)=w$),  then inner variation on the ``$h$ side'',  gives us the existence of a {\em holomorphic Ahlfors-Hopf differential},
\begin{equation}\label{AHDifferential}
\Phi=\mA'(\IK(w,h))h_w\overline{h_\wbar}\,\eta(h).
\end{equation}
We remark that the existence of the minimiser does depend on the regularity of $\mA$. For instance if the functions are in the exponential class, that is there exists constants $c_0>0$ and $p>0$ such that
\begin{equation}\label{1.5}
\mA(t)\geq c_0\, e^{pt},
\end{equation}
then there is a homeomorphic minimiser \cite{MY5}. However, if $\mA$ satisfies
\begin{equation}\label{1.6}
\mA(t)\geq c_0\, t^p,
\end{equation}
where $c_0>0$ and $p>1$, then a minimiser exists and  lives in the enlarged space $\mF_{\mA,\eta}$, which we define next.

\subsection{The enlarged spaces $\mathsf{F}_{\mA,\eta}$ and $\mathsf{H}_{\mA,\eta}$.}  Let $f_0:\oD\to\oD$ be a homeomorphism with $\mathsf{E}_{\mA,\eta}(f_0)<\infty$, and set $h_0=f_0^{-1}$. We then define $\mathsf{F}_{\mA,\eta}$ as the space of functions $f:\ID\to\ID$ satisfying the following three conditions:
\begin{enumerate} 
\item $f\in W^{1,q}(\ID,\ID)$ for some $q=\frac{2p}{p+1}$, where $p$ is from (\ref{1.6}), and
\[
\|Df\|_{L^q(\ID)}\leq\|Df_0\|_{L^q(\ID)}+1,
\]
\item There exists $h\in C(\overline{\ID})$, $h|_{\partial\ID}=h_0|_{\partial\ID}$, $h$ has finite distortion, and
\[
\mathsf{E}_{\mA,\eta}^\ast(h)
\leq\int_\ID\mA(\IK(w,h_0))J(w,h_0)\eta(h_0)dw+1.
\]
\item There is a dense set $X\subset\ID$ of full measure, $|\ID-X|=0$, such that  
\[ h\circ f(z)=z, \quad \mbox{ for every $z\in X$} \]
and $J(w,h)=0$ for almost every $w\in\ID-f(X)$.
\end{enumerate}

We then define $\mH_{\mA,\eta}$ to be the space of all such pseudo-inverses $h$.   

\subsection{Smoothness and Uniqueness.} 
The proof for the existence of a minimiser in $\mF_{\mA,\eta}$ is essentially same as our earlier work for the particular case $\mA(t)=t^p$ in \cite[\S 5]{MY1}. Moreover, if 
\begin{equation}\label{1.7}
t\mA'(t)\leq c_1\mA(t),
\end{equation}
then the existence of the Ahlfors-Hopf differential is guaranteed, \cite{IMO}.
\begin{theorem}
If $\mA$ satisfies (\ref{1.6}), then there is a minimiser in the enlarged space $\mF_\mA$. Moreover, if (\ref{1.7}) is also satisfied, then there exists a holomorphic Ahlfors-Hopf differential $\Phi\in L^1(\ID)$ which satisfies (\ref{AHDifferential}).
\end{theorem}

Our previous results concerning diffeomorphic minimisers (see \cite{MY1},\cite{MY3}) were mainly based on an analysis of the $f$-side equation (\ref{1.3}). However, to achieve them,   higher regularity assumptions were required. We later developed new methods based on the $h$-side Ahlfors-Hopf differential, and proved that the minimisers of exponential distortion between analytically finite surfaces must be diffeomorphisms, \cite{MY5}. The following result follows directly from the same considerations as in \cite[Theorem 18]{MY5}:
\begin{theorem}\label{ExistenceDiff} (Existence of diffeomorphic solutions)
Let $\Phi\in L^1(\ID)$, and suppose there exits a constant $c_0>0$ such that
\begin{equation}\label{1.8}
c_0\mA(t)\leq t\mA'(t),
\end{equation}
Then there exists a $g\in\mF_{\mA,\eta}$ such that its pseudo-inverse $H$ satisfies (\ref{AHDifferential}) in $\ID$, and that $g$ is locally diffeomorphic in $\ID$.
\end{theorem}
It is to be noted that we have not claimed that $g$ is surjective.

For functions in the exponential class, we can say even more, see \cite{MY5}:
\begin{theorem}
If $\mA$ satisfies (\ref{1.5}), and the minimiser $f$ satisfies (\ref{AHDifferential}), then $f$ is diffeomorphic in $\ID$.
\end{theorem}

We remark the assumption that the minimiser $f$ satisifies (\ref{AHDifferential}) is not automatically true, \cite{MY1}. In fact, in the surface case \cite{MY5}, we applied the Riemann-Roch theorem to prove the existence of the Ahlfors-Hopf equation. The key to travel from Theorem 2 to Theorem 3 is to prove a uniqueness theorem to get $h\equiv H$, where $H=g^{-1}$, $h=f^{-1}$. 

Here we prove the following uniqueness theorem in the unweighted case given that one solution is regular:
\begin{theorem}\label{UniqueDiff}
Let $h$ and $H$ both be in the enlarged space $\mH_\mA$, where the boundary data $h|_\IS$ and $H|_\IS$ are both monotone but might be different, $h$ and $H$ both satisfy the same Ahlfors-Hopf differential (\ref{AHDifferential}), where $\eta\equiv1$, and $h(w_0)=H(w_0)=z_0$, $h(w_1)=H(w_1)=z_1$, where $w_0,z_0\in\ID$, $w_1,z_1\in\IS$. Suppose further that the pseudo-inverse $g$ of $H$ is diffeomorphic in $\ID$. Then $h\equiv H$ in $\oD$.
\end{theorem}
We remark that in the above theorem we have extended the definition of $\mH_\mA$ to allow topologically monotone boundary values.  We next prove the uniqueness of homeomorphic minimisers.

\begin{theorem} \label{Uniqueness} (Uniqueness of homeomorphic solutions) Fix homeomorphic boundary data $f_0$, and let $\mA$ satisfy (\ref{1.7}).
Let $h\in\mH_\mA$ be a homeomorphic minimiser for the problem (\ref{1.2}). If $H\in\mH_\mA$ satisfies the holomorphic Ahlfors-Hopf equation
\begin{equation}
\Phi=\mA'(H)\,H_w\,\overline{H_\wbar}\;\eta(H)\in L^1(\ID),
\end{equation}
then $h\equiv H$ in $\oD$.
\end{theorem}

Putting the above results together, we can now give the following result for the unweighted problem (\ref{1.1}). Note this improves our earlier results in \cite{MY1} and \cite{MY3}.

\begin{theorem}
Consider the unweighted problem (\ref{1.1}). Assume that $\mA$ satisfies (\ref{1.7}) and (\ref{1.8}). Then there exists a minimiser $f:\ID\to f(\ID)$ which  is a diffeomorphism. Furthermore, if $h$ is a self-homeomorphism of $\ID$, then $h$ is the unique minimiser in $\mH_\mA$.
\end{theorem}

We now consider the case $\mA$ is $L^p$-like (i.e. (\ref{1.7})-(\ref{1.8}) hold), but $\eta$ is not necessarily a constant. For the case that a minimiser $f\in\mF_\mA$ is not a self-homeomorphism of $\oD$, we can also say something. 

Note that in this case, its pseudo-inverse $h$ is a topologically monotone mapping. For each $z\in\oD$, $h^{-1}(z)$ is a closed connected set; since $h$ is diffeomorphic from $\Omega=f(\ID)$ onto $\ID$, these bad points appear only on $\IS$. Now if for a point $z\in\IS$, $h^{-1}(z)$ is not a singleton, then it must be a closed connected set such that $h^{-1}(z)\cap\IS$ is a single point. We call this set $h^{-1}(z)$ a {\em hair} and $h^{-1}(z)\cap\IS$ the tip of the hair. The set $E=\ID\setminus\Omega$ is then a union of these hairs, and we call the set $Z=\{z\in\IS:h^{-1} \mbox{ is not a singleton}\}$. By the definition of $\mH_{\mA,\eta}$, $J(w,h)=0$ almost everywhere in $E$. Also, $h$ satisfies the holomorphic Ahlfors-Hopf differential (\ref{AHDifferential}). Note then $J(w,h)\neq0$ almost everywhere in $\ID$. Thus $|E|=0$. We prove the following:

\begin{theorem}\label{hairs}
Let the notation be as stated above. Then the set $Z$ has Hausdorff measure  $\mathcal{H}^1(Z)=0$.
\end{theorem}

\begin{center}
\scalebox{0.25}{\includegraphics*{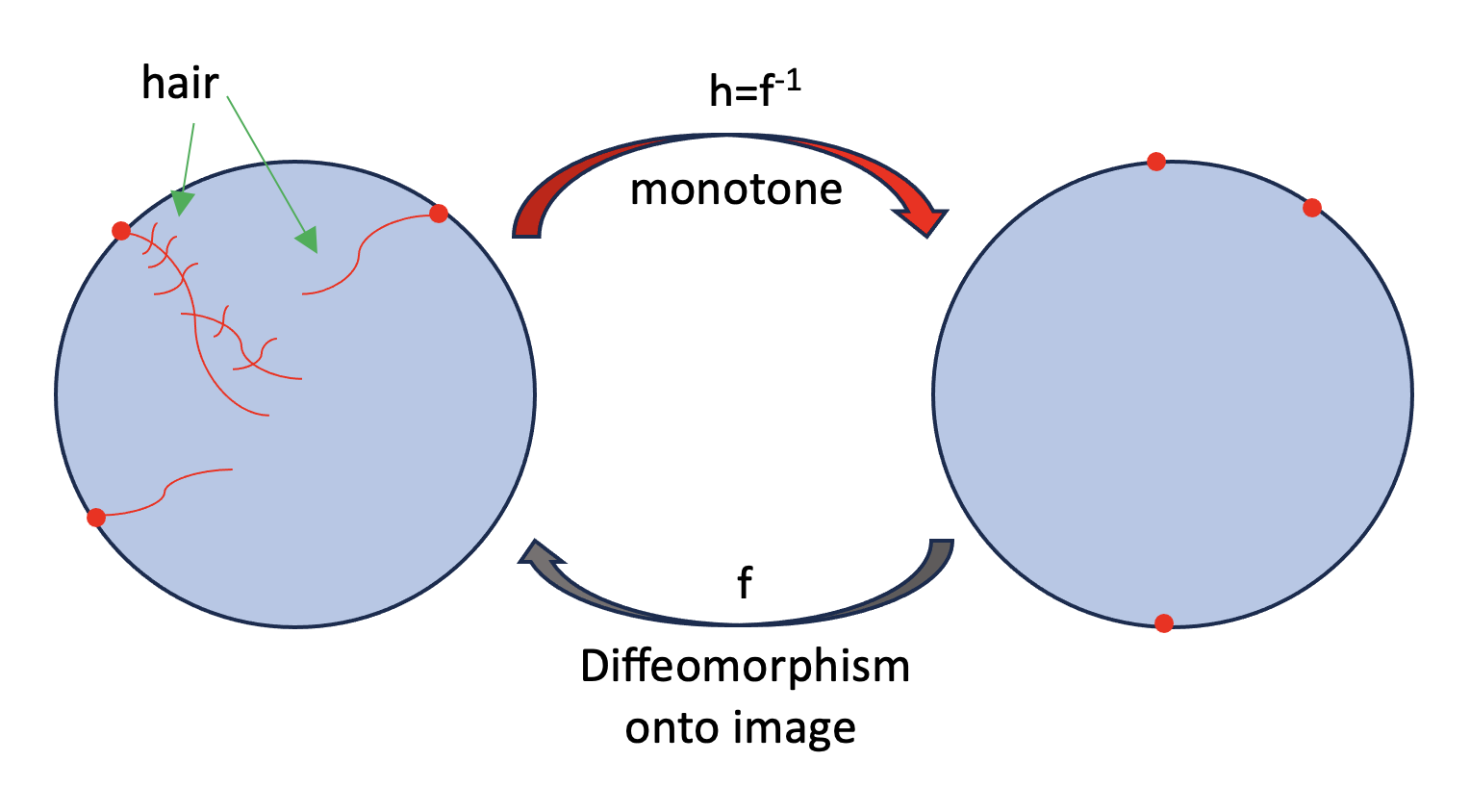}}
\end{center}
\noindent{\bf Figure 1.} {\it  A typical monotone minimiser $h$ which shrinks hairs.  As $h$ is a uniform limit of homeomorphisms,  each hair is cellular;  its complement in $\oC$ is both connected and simply connected.}

\medskip

In Section \ref{MonotoneExample}, we give an example of a weighted case which is induced by a finite distortion problem on surfaces, where the uniqueness applies and the minimiser is a diffeomorphism.

\section{Diffeomorphism; proof of Theorem \ref{UniqueDiff}}

From our hypotheses we have a minimiser $f$ to the boundary value problem (\ref{1.1}) with monotone pseudo-inverse $h$ which satisfies the Ahlfors-Hopf equation (\ref{AHDifferential}) with $\eta\equiv1$,  along with a solution $H$ to the same Ahlfors-Hopf equation as $h$, and $H$ satisfies the conditions as in Theorem \ref{UniqueDiff}. Moreover, they have a common pair of points with distinguished values, for notational ease we write $h(0)=H(0)=0$, and $h(1)=H(1)=1$. Our task of this section is to prove $h=H$ in $\oD$.\\

Let $F=h-H$, and $\Omega\subset\ID$ be the set where $H$ is diffeomorphic. The following lemma follows from a similar argument as given in \cite[Lemma 6]{MY5}.  The proof is simply to write the Ahlfors-Hopf equation as a degenerate elliptic Beltrami equation $h_\zbar={\cal B}(h_z)$, show that ${\cal B}$ satisfies local elliptic bounds away from $0$ from which it follows with a few estimates that $F$ is locally the solution to an elliptic Beltrami equation.  
\begin{lemma}
$F$ is locally quasiregular in $\Omega$. In particular, $F$ is open in $\Omega$.
\end{lemma}
Further note that the growth estimates on $\mA$ and as $\mathsf{E}_{\mA,\eta}^\ast(h)$ and $\mathsf{E}_{\mA,\eta}^\ast(H)$ are finite,  we have $h,H\in W^{1,2}(\ID)$.

\medskip

Now since $F=h-H\in W^{1,2}(\ID)$ and $|\Omega^c|=0$, we find that by \cite{GV}:
\begin{lemma}
$F$ is a $W^{1,2}(\ID)$ finite distortion function. In particular, $F$ satisfies Lusin's $\mathcal{N}$ condition.
\end{lemma}

We now start to inspect the difference on the boundary. Recall both $h$ and $H$ map the unit circle $\IS$ to itself, and they are both monotone ($h$ is in fact homeomorphic on $\IS$). Consider two functions $\alpha(\theta):[0,2\pi)\to[0,2\pi)$ and $\beta(\theta):[0,2\pi)\to[0,2\pi)$, where $\alpha$ and $\beta$ are both sense-preserving monotone mappings onto the interval $[0,2\pi)$. We write the difference $F(e^{i\theta})=e^{i\alpha(\theta)}-e^{i\beta(\theta)}$ and consider the set $S=\{\theta:\alpha(\theta)=\beta(\theta)\}$. Note $0\in S$ so it is at least non-empty. From the continuity of $\alpha-\beta$, $S$ is a closed set and this means the complement $[0,2\pi)\setminus S$ is a countable union of non-intersecting open intervals. Let $(\theta_1,\theta_2)$ be such an interval. The image of $e^{i\alpha(\theta)}-e^{i\beta(\theta)}$ is a ``lobe'' as on each interval the increment of the argument is monotonically increasing or decreasing. Let $\gamma$ be the image of $e^{i\alpha(\theta)}-e^{i\beta(\theta)}$, $\theta\in(\theta_1,\theta_2)$. Now $\gamma\cup\{0\}$ is a Jordan curve, and $\gamma^\circ$ (the region bounded by $\gamma$) is a simply-connected open set with $0$ on its boundary; a lobe $L$, and its closure is $\overline{L}=L\cup\gamma\cup\{0\}$.
\begin{lemma}
Let $\gamma=F(e^{i\theta})=e^{i\alpha(\theta)}-e^{i\beta(\theta)}$, $\theta\in(\theta_1,\theta_2)$ define a lobe, i.e. $\alpha(\theta_1)=\beta(\theta_1)=\tilde{\theta}_1$, $\alpha(\theta_2)=\beta(\theta_2)=\tilde{\theta}_2$ and either $\alpha>\beta$ throughout $(\theta_1,\theta_2)$ or $\alpha<\beta$ throughout $(\theta_1,\theta_2)$. Then along $\gamma$ we have
\begin{enumerate}[label=(\roman*)]
\item The argument $\arg(F(e^{i\theta}))$ is monotone.
\item For each point $re^{i\delta}=F(e^{i\theta})$, the line segment $\{\rho e^{i\delta}:r\in(0,\rho)\}\subset L$.
\item The variation of argument $(\varg)$ of $\gamma$ is $\tilde{\theta}_2-\tilde{\theta}_1$.
\end{enumerate}

\end{lemma}
\noindent{\bf Proof.} We assume $\alpha$, $\beta$ are $C^1$-smooth, then
\begin{align*}
\frac{d}{d\theta}\tan\left(\arg(F(e^{i\theta}))\right)&=\frac{d}{d\theta}\left(\frac{\sin\alpha(\theta)-\sin\beta(\theta)}{\cos\alpha(\theta)-\cos\beta(\theta)}\right)\\
&=\frac{(\alpha'+\beta')(1-\cos(\alpha-\beta))}{(\cos\alpha-\cos\beta)^2}\geq0.
\end{align*}
So it is non-decreasing, which proves (i). (ii) is a direct corollary of (i) and the fact $0\in\partial L$. To prove (iii), by (i), we only need to consider the tangent lines at $\theta_1$ and $\theta_2$, i.e.
\[
\varg(e^{i\alpha}-e^{i\beta})=\arg[(e^{i\alpha}-e^{i\beta})'(\theta_2)]-\arg[(e^{i\alpha}-e^{i\beta})'(\theta_1)].
\]
In fact we can compute them directly:
\[
(e^{i\alpha}-e^{i\beta})'(\theta_1)=ie^{i\alpha(\theta_1)}\alpha'(\theta_1)-ie^{i\beta(\theta_1)}\beta'(\theta_1)=ie^{i\tilde{\theta}_1}(\alpha'(\theta_1)-\beta'(\theta_1)),
\]
\[
(e^{i\alpha}-e^{i\beta})'(\theta_2)=ie^{i\alpha(\theta_2)}\alpha'(\theta_2)-ie^{i\beta(\theta_2)}\beta'(\theta_2)=ie^{i\tilde{\theta}_2}(\alpha'(\theta_2)-\beta'(\theta_2)).
\]
If either $\alpha$ or $\beta$ is not smooth, we may approximate it by smooth functions $\alpha_\epsilon\to\alpha$, or $\beta_\epsilon\to\beta$, which is uniform in $C^0$, and the images of the end points are kept. By uniform convergence the condition $\alpha>\beta$ (or $\beta>\alpha$) is kept, and all the results (i), (ii), (iii)  are also kept.\hfill$\Box$

\bigskip

\begin{center}
\scalebox{0.25}{\includegraphics*{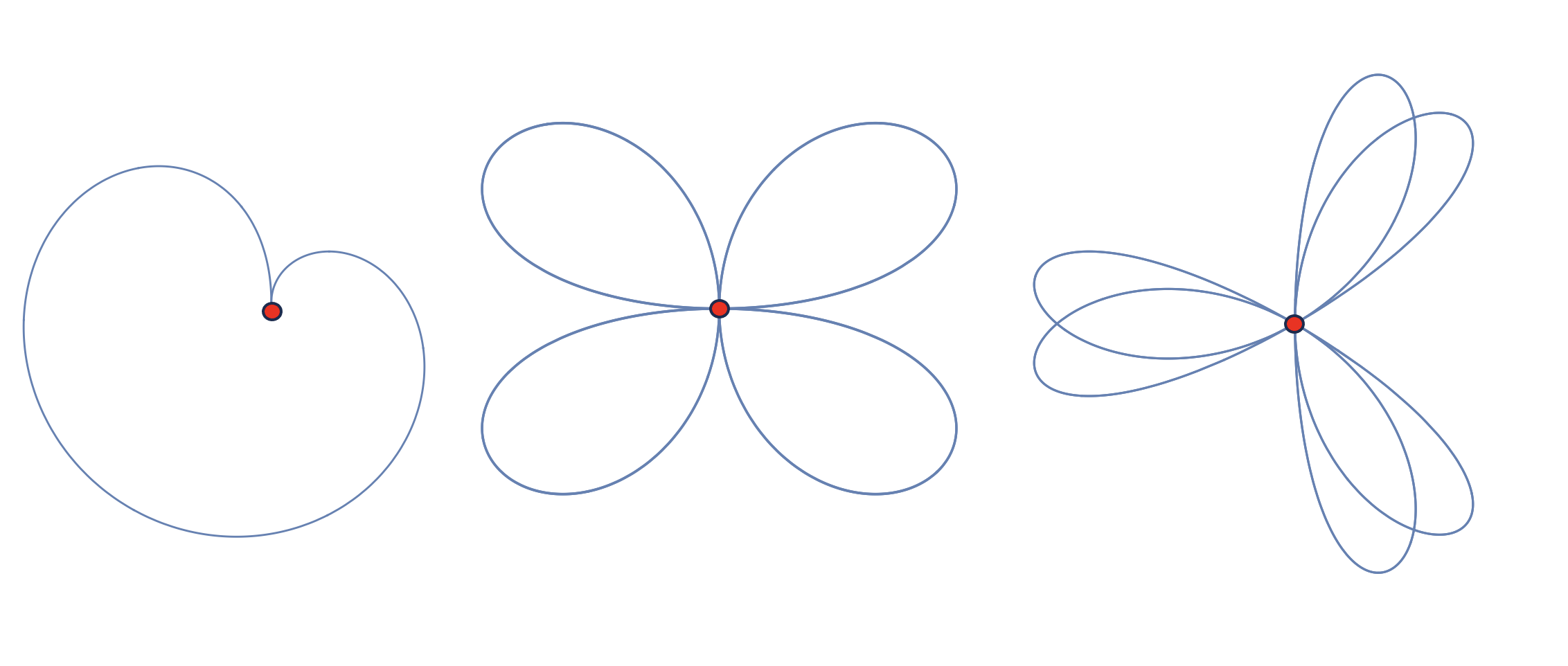}} \\
\end{center}
\noindent{\bf Figure 2.} {\it The possible behaviours of the collection of lobes. From left: one lobe, four lobes, six overlapping lobes.}

\medskip

From now we enumerate the lobes $L_i$ and let $L=\bigcup L_i$ be the collection of all these lobes. 

\begin{lemma} With the notation above, 
\begin{enumerate}[label=(\roman*)]
\item The total variation of argument for $F(\IS)$ satisfies $\varg(F)\leq2\pi$.
\item $0\in\partial L$.
\item For each $\theta\in[0,2\pi)$, there is an $r\in[0,\infty)$ such that 
\[ \{s\in[0,\infty]:se^{i\theta}\in L\}=(0,r).\]Furthermore,  $r(\theta)e^{i\theta}:\theta\in[0,2\pi)$ is a continuous curve.
\item If a bounded domain $A$ has boundary $\partial A\subset F(\IS)$, then $\overline{A}\subset\overline{L}$
\end{enumerate}
\end{lemma}

\noindent{\bf Proof.} For (i), by (iii) of Lemma 1, we have $\varg(F)=2\pi-|X|$, where
\[
X=\{\tilde{\theta}\in[0,2\pi]:\alpha(\theta)=\beta(\theta)=\tilde{\theta}\}.
\]
(ii) is a direct result of (i); (iii) follows from (ii) of lemma 3 and the continuity of $F$. For (iv), we have the cone $C:=\{\rho e^{\theta}:\theta\in(\min\arg(A),\max\arg(A)),\rho\in(0,r(\theta))\}$ such that $A\subset C\subset L$.\hfill$\Box$

\medskip

We remark the last result could be even better. Consider $\overline{A}=A\cup\partial A$, a compact set. Its complement in the extended plane, $\oC\setminus\overline{A}$ is composed of countably many open domains. There exists a unique component of it that contains the point infinity -- call it $\overline{A}^c_\infty$.
\begin{lemma}
If $\partial\overline{A}^c_\infty\subset F(\IS)$, then $\overline{A}\subset \overline{L}$.
\end{lemma}
\noindent{\bf Proof.} We apply (iv) of Lemma 4 to $(\overline{A}^c_\infty)^c$. Note then $\overline{A}\subset\overline{(\overline{A}^c_\infty)^c}\subset\overline{L}$.\hfill$\Box$

\medskip

We now consider the hairs inside $\ID$. As we discussed above, each hair $E_i$, ($i\in\mathcal{I}$) starts from a point $\xi_i$ on $\IS$, and $H(E_i)=H(\xi_i)$ is a constant. Note $\IS\cup\left(\bigcup_{i\in\mathcal{I}}E_i\right)=\oD-\Omega$ is a closed set with measure $0$. In particular, it is nowhere dense in $\IC$. However, there might be uncountably many hairs, and $\{\xi_i\}_{i\in\mathcal{I}}$ could be an uncountable set on $\IS$.\\

In the following lemma, we use $F(\IS)^c_\infty$, $F(\oD)^c_\infty$ to denote the component of the complement of $F(\IS)$ and $F(\oD)$ in $\oC$ which contains the point  $\infty$, respectively.
\begin{lemma}
If there is a point $z\in\ID$ such that $F(z)\in\partial F(\oD)^c_\infty$, then $F(z)\in F(\IS)$.
\end{lemma}
\noindent{\bf Proof.} Assume there is such a point $z\in\ID$ so that $F(z)\in\partial F(\overline{\ID})^c_\infty\setminus F(\IS)$. Then $F^{-1}(F(\IS)^c_\infty)\subset\ID$ is an open neighbourhood of $z$, and the winding number $w(F,F(\IS)^c_\infty)$ is a constant in this set. Precisely, the winding number is defined by
\[
w(F,F(\IS)^c_\infty)=\deg(F,x)=\sum_{a\in F^{-1}(x)}\mbox{sign}(J(a,F))
\]
which is independent of the choice of $x\in F(\IS)^c_\infty$. See \cite[p.156, Proposition 4.2]{OR}. However, on one hand, we have $\deg(F,\infty)=0$; on the other hand, by the Lusin $\mathcal{N}$ property, there must be such a non-critical value $x\in F(\ID)^c_\infty$ so that $J(a,F)>0$ for all $a\in F^{-1}(x)$ ($a\in\Omega$ and $F$ is open at $a$), so $\deg(F,x)\geq1$. Thus there cannot be such a point $z$.\hfill$\Box$

\bigskip

\begin{center}
\scalebox{0.4}{\includegraphics*{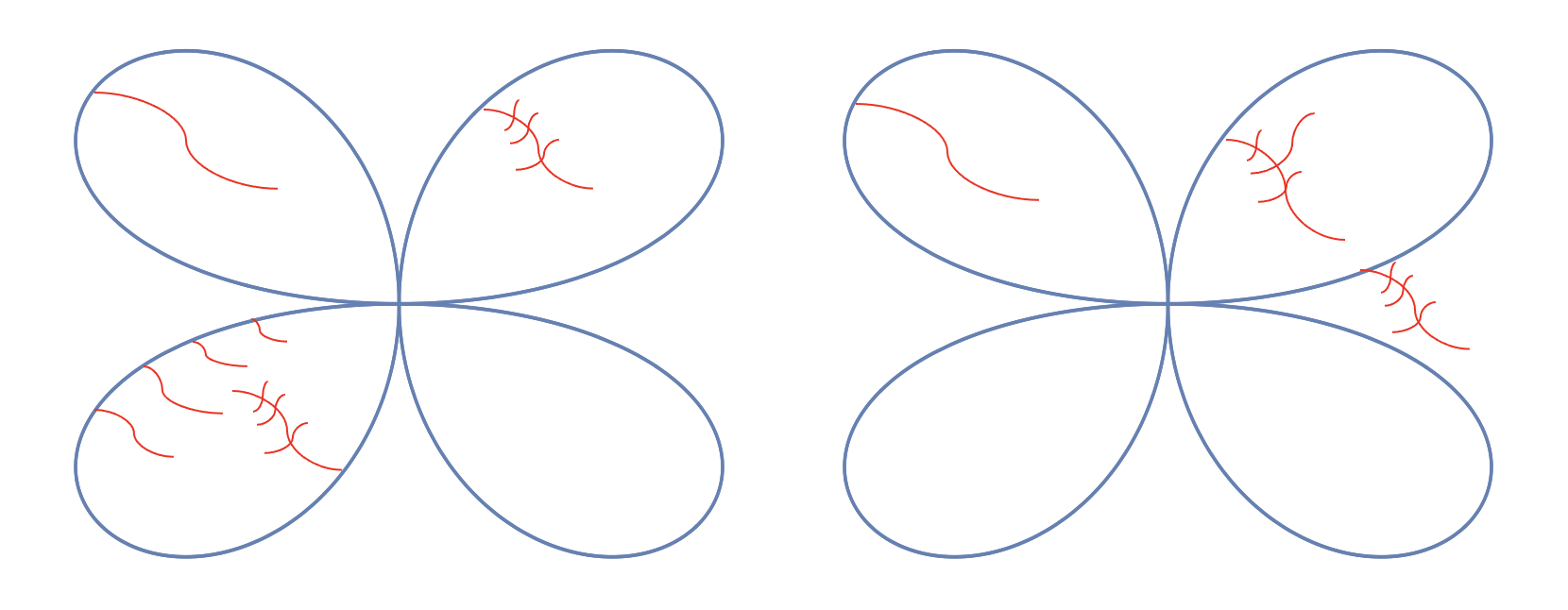}} \\
\end{center}
\noindent{\bf Figure 3.} {\it  The image of the hairs cannot lie outside the lobes $L$: Left OK, Right not possible.}
\medskip

\noindent{\bf Proof of Theorem \ref{UniqueDiff}.} Note $h(0)=H(0)=0$, so $0\in\Omega$ is not on any hair. In particular, if $F$ is not a constant in $\Omega$, it is open in a neighbourhood of $0$. However, by Lemma 5 and Lemma 6, $F(\oD)\subset\overline{L}$, thus $F(0)=0\in\partial F(\oD)$, which is then impossible. Thus $F\equiv0$ in $\oD$.\hfill$\Box$

\medskip

\section{Uniqueness; proof of Theorem \ref{Uniqueness}} 
We now assume a $L^p$ minimiser $h:\oD\to\oD$ is homeomorphic. By what we have established above, it must in fact be a diffeomorphism. In this section we shall prove it must be unique. In fact we will prove   a slightly more general result.
\begin{theorem}
Let $h\in\mH_p$ be a minimiser for the $L^p$ problem (\ref{1.1}) with an associated homeomorphic minimising sequence.  That is suppose there exists a sequence $h_n\in\mH_p$, $n\geq 1$ and each $h_n:\oD\to\oD$ is a homeomorphism and  that
\[
\mathsf{E}_p^\ast(h_n)\to\mathsf{E}_p^\ast(h)=\min_{g\in\mathsf{H}_p}\mathsf{E}_p^\ast(g).
\] 
Then $h$ is a unique minimiser in $\mathsf{H}_p$. In fact, $h$ is the unique critical point for the inner variational equation, that is it is the only function in $\mathsf{H}_p$ that satisfies the equation
\begin{equation}\label{4.1}
\int_\ID\IK_h^{p-1}\,h_w\,\overline{h_\wbar}\, \eta(h)\, \varphi_\wbar=0,\quad\forall\varphi\in C_0^\infty(\ID).
\end{equation}
\end{theorem}
\noindent{\bf Proof.} Assume that $g\in\mathsf{H}_p$ satisfy (\ref{4.1}). Then it satisfies the Ahlfors-Hopf equation (\ref{1.2}) for some holomorphic function $\Phi$. By \cite{MY4}, for each $h_n$, $\mathsf{E}_p^\ast(g)\leq\mathsf{E}_p^\ast(h_n)$, thus $\mathsf{E}_p^\ast(g)\leq\lim_{n\to\infty}\mathsf{E}_p^\ast(h_n)=\mathsf{E}_p^\ast(h)$. However, $h$ is a minimiser, so $\mathsf{E}_p^\ast(g)=\mathsf{E}_p^\ast(h)$. Again by \cite{MY4}, this happens only when $g=h$.\hfill$\Box$

\medskip

We remark that in \cite{MY4}, the proof depends on the Reich-Strebel inequality
\begin{equation}\label{RSinequality}
\int_\ID|\Phi|\leq\int_\ID\sqrt{|\Phi(\xi_n)|}\sqrt{|\Phi|}\left|(\xi_n)_z-\frac{\Phi}{|\Phi|}(\xi_n)_\zbar\right|,
\end{equation}
where $\xi_n=h_n^{-1}\circ g$. To get (\ref{RSinequality}), the functions $h_n$ and $g$ were both assumed to be homeomorphisms. However, once we inspect the proof of the Reich-Strebel inequality, we find that $g$ being in the extended space $\mathsf{H}_p$ is allowed. See the illustration below.

\begin{center}
\scalebox{0.12}{\includegraphics*{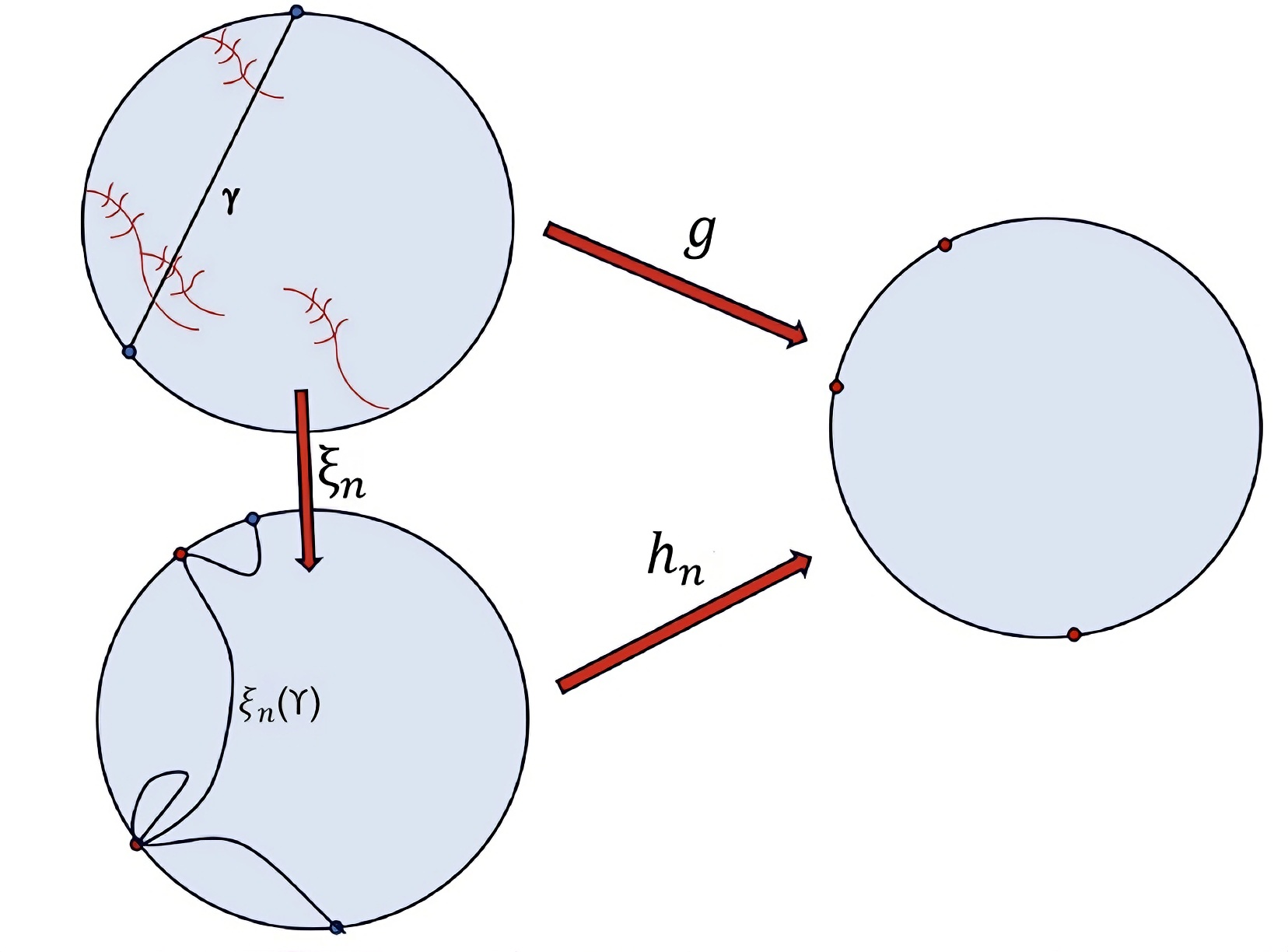}} \\
\end{center}
\noindent{\bf Figure 4.} {\it If $g\in\mathsf{H}_p$ is not a homeomorphism, we still have $|\xi_n(\gamma)|\geq|\gamma|$ for any $|\Phi|^{1/2}$-geodesic $\gamma$.}

\medskip

See \cite[\S 2.2]{IO} for how this produces the Reich-Strebel inequality.
\medskip

\section{The hairs; proof of Theorem \ref{hairs}}
Recall that $E$ has measure $|E|=0$. Now for almost every circle $S_r=S(0,r)$ inside $\ID$, $\mathcal{H}^1(S_r\cap E)=0$, thus we can choose  a (countable) sequence $r_n\to1$ so that $|S_{r_n}\cap E=0|$ for each $n$. Note then each hair must hit some $S_{r_n}$ for all sufficiently large $n$, so $Z=\bigcup_{n=1}^\infty Z_n$, where $Z_n:=h(S_{r_n}\cap E)$. On each $S_{r_n}$, since $h$ is Lipschitz \cite{IKO}, $|Z_n|\leq|S_{r_n}\cap E|=0$. Thus
\[
|Z|=|\bigcup_{n=1}^\infty Z_n|=0.
\]
This proves Theorem \ref{hairs}.

\section{Monotone automorphic mappings.}\label{MonotoneExample}

There is a natural example where the hairs described above cannot occur,  that is they must not run out to the boundary.  This is the case of lifts of monotone maps between closed Riemann surfaces.  Let $\Gamma_1$ and $\Gamma_2$ be compact torsion Fuchsian groups acting on the hyperbolic plane.  Thus the orbit spaces $\Sigma_1=\IH^2/\Gamma_1$ and $\Sigma_2=\IH^2/\Gamma_2$ are closed Riemann surfaces of genus at least two.  If they have the same genus then these surfaces are diffeomorphic.  This diffeomorphism acts as a barrier for a minimisation problem we subsequently discuss.  

First we will present a few lemmas we need,  this all closely follows Ahlfors \cite[\S V, (30) and Theorem 13]{Ah3}.

\begin{lemma}\label{Projection} Let $\Gamma$ be a Fuchsian surface group and $\tilde{\rho}:\ID\to\IR$ a continuous  automorphic function
\[ \tilde{\rho}(\gamma(z))=\tilde{\rho}(z).\]
Then $\tilde{\rho}$ projects to a well defined continuous function $\rho:\Sigma=\ID/\Gamma\to \IR$ and 
\begin{equation}
\int_\Sigma  \rho (w)\; d\sigma(w) = \int_{P} \tilde{\rho}(z)\; \frac{dz}{(1-|z|^2)^2},
\end{equation}
where $\sigma$ is the hyperbolic area measure on $\Sigma$ and $P$ is any measurable fundamental region.
\end{lemma}
\noindent {\bf Proof.} That $\tilde{\rho}$ induces $\rho$ and is well defined is simply covering space theory. The formula follows since it clearly holds for locally constant simple functions as the projection is an isometry of the two hyperbolic metrics (and measures),  and these are uniformly dense in the continuous functions. \hfill $\Box$

\medskip
\begin{lemma}\label{tAlpha} The function
\begin{equation} 
\tilde{\alpha}(z) = \sum_{\gamma\in\Gamma_1} |\gamma'(z)|^2\; (1-|z|^2)^2 =\sum_{\gamma\in\Gamma_1}  (1-|\gamma(z)|^2)^2
\end{equation}
is continuous and automorphic.
\end{lemma}
\noindent{\bf Proof.} The uniform convergence of the Poincar\'e series $\sum_{\gamma\in\Gamma_1} |\gamma'(z)|^2$ for Fuchsian groups of the first kind is well known.  So continuity follows.  We also compute that for $\beta\in\Gamma_1$
\begin{eqnarray*}
\tilde{\alpha}(\beta(z)) &=& \sum_{\gamma\in\Gamma_1} |\gamma'(\beta(z))|^2\; (1-|\beta(z)|^2)^2 \\
&=& \sum_{\gamma\in\Gamma_1} |\gamma'(\beta(z))|^2|\beta'(z)|^2\; \frac{(1-|\beta(z)|^2)^2}{|\beta'(z)|^2} \\
&=& \sum_{\gamma\in\Gamma_1} |(\gamma\circ\beta)'(z))|^2 \;  (1-|z|^2)^2  \\
&=& \sum_{\gamma\in\Gamma_1} |(\gamma)'(z))|^2 \;  (1-|z|^2)^2  = \tilde{\alpha}(z).
\end{eqnarray*}
The last line following because the sums are the same. \hfill $\Box$

\medskip

By Lemma \ref{Projection} and Lemma \ref{tAlpha}, there is a function $\alpha:\Sigma_1=\ID/\Gamma_1\to\IR$ which lifts to $\tilde{\alpha}:\ID\to\IR$. We next show this also works with a multiplication of $\IK^p(z,f)$.

\begin{lemma} Suppose that $f:\Sigma_1\to\Sigma_2$ is a diffeormorphism of finite distortion inducing an isomorphism of fundamental groups,
and let $\tilde{f}$ be the lift of $f$ to the universal cover $\ID$,
\[ 
\tilde{f}\circ \Gamma_1 = \Gamma_2 \circ \tilde{f}. 
\] Then
\begin{equation}
\IK^p(z,\tilde{f})\tilde{\alpha}(z)
\end{equation}
is a smooth automorphic function.
\end{lemma} 
\noindent{\bf Proof.} With the obvious notation we calculate
\begin{eqnarray*}
\IK^p(\beta_1(z),\tilde{f})  \tilde{\alpha}(\beta_1(z)) &=&\IK^p(z,\widetilde{f\circ \beta_1})  \; \tilde{\alpha}(z) = \IK^p(z,\tilde{f}\circ \beta_1) \; \tilde{\alpha}(z)  \\
&=& \IK^p(z,\beta_2\circ \tilde{f} ) \; \tilde{\alpha}(z) = \IK^p(z,\tilde{f} ) \; \tilde{\alpha}(z).
\end{eqnarray*}
This is what we wanted to prove.\hfill $\Box$

\medskip

\begin{theorem}\label{thm9} For a diffeomorphism $f:\Sigma_1\to\Sigma_2$,
\begin{equation}\label{5.4}
\int_{\Sigma_1} \IK^p(z,f) \alpha(z) \;d\sigma(z) = \int_\ID \IK^p(z,\tilde{f}) \; dz.
\end{equation}
Further,  with $h=f^{-1}$,  we have 
\begin{equation}\label{5.5}
\int_{\Sigma_2} \IK^p(w,h)\;  J(w,h) \; \alpha(h) \;d\sigma(w) = \int_\ID \IK^p(w,\tilde{h})\, J(w,\tilde{h}) \; dw.
\end{equation}
\end{theorem}
\noindent{\bf Proof.} We prove (\ref{5.4}) as the second part follows by the change of variable formula.  We calculate that for any fundamental polygon $P$
\begin{eqnarray*} 
\lefteqn{\int_{\Sigma_1} \IK^p(z,f) \alpha(z) \;d\sigma(z)}\\ 
& = &  \int_P \IK^p(z,\tilde{f})  \tilde{\alpha}(z) \;\frac{dz}{(1-|z|^2)^2} \\
& = &  \int_P \IK^p(z,\tilde{f}) \sum_{\gamma_1\in\Gamma} |\gamma_1'(z)|^2\; (1-|z|^2)^2  \; \frac{dz}{(1-|z|^2)^2} \\
& = &  \int_P \IK^p(z,\tilde{f}) \sum_{\gamma_1\in\Gamma} |\gamma_1'(z)|^2  \;dz  = \sum_{\gamma_1\in\Gamma}  \int_P \IK^p(z,\tilde{f}) |\gamma_1'(z)|^2  \;dz  \\
& = & \sum_{\gamma_1\in\Gamma}  \int_{\gamma_1(P)} \IK^p(\gamma_1^{-1},\tilde{f})  \;dz  =  \int_{\ID} \IK^p(\gamma_1^{-1},\tilde{f})  \;dz  =  \int_{\ID} \IK^p(z,\tilde{f}\circ\gamma_1^{-1})  \;dz  \\
& = &   \int_{\ID} \IK^p(z,\gamma_2^{-1}\circ \tilde{f})\;dz  =  \int_{\ID} \IK^p(z,\tilde{f})  \;dz.
\end{eqnarray*}
This completes the proof. \hfill $\Box$

\medskip

With Theorem \ref{thm9} a question about extremal mappings between Riemann surfaces now reduces to our previous cases about extremal mappings between the disk.  However there is a key difference.  Instead of minimising over all mappings $f:\ID\to\ID$ with the same boundary values,  we must minimise over all diffeomorphisms automorphic with respect to the associated Fuchsian groups.  This too is discussed in \cite{Ah3} and one obtains the existence of an automorphic holomorphic quadratic differential $\Phi:\ID\to \ID $,
\begin{eqnarray}\Phi(w) &=& \IK^{p-1}(w,\tilde{h}) \, \tilde{h}_w(w) \, \overline{\tilde{h}_\wbar},\\
\label{Phi} \Phi(w)&=&\Phi(\gamma_2(w))(\gamma_2'(w))^2, \quad \gamma_2\in \Gamma_2. \end{eqnarray}
The idea is simply to obtain the same variational equations by using a local inner variation supported inside any fundamental polygon - the automorphy condition in fact guarantees the same boundary values if normalised, and otherwise up to a M\"obius element of the group $\Gamma_1$ (or $\Gamma_2$ for $\tilde{h}$).  

Another more serious problem is that for the global problem the holomorphic quadratic differential lies in $L^1(\ID)$,  however the automorphic constraint (\ref{Phi}) guarantees $\Phi\not\in L^1(\ID)$ unless $\Phi\equiv 0$ and $\tilde{h}$ is conformal.  To see this we calculate that 
\begin{eqnarray*}
 \int_\ID |\Phi(w)|\; dw &=&\int_\ID |\Phi(\gamma_2(w))| |\gamma_2'(w)|^2 \\
 &=&\sum_{\gamma_2}\int_{\gamma_2(P)} |\Phi(\gamma_2(w))| |\gamma_2'(w)|^2 \\
 &=&\sum_{\gamma_2}\int_{P} |\Phi(w)| \; dw =+\infty.
\end{eqnarray*}

The same modulus of continuity estimates and the fact that the automorphy condition is preserved under local uniform limits implies that a minimising sequence of diffeomorphisms for the functional (\ref{5.5}) has a monotone limit. Ahlfors actually established these modulus of continuity estimates first,  at least for $p>2$ and in a less refined form.  The arguments we discussed earlier now give minimisers subject to being  automorphic with respect to the groups $\Gamma_1$ and $\Gamma_2$.  We summarise this as follows.

\begin{lemma}\label{lem10} Let $\{h_j:\Sigma_2\to\Sigma_1\}$ be a minimising sequence of diffeomorphisms for the problem 
\begin{equation}
\inf_{h \; {\rm diffeo}} \int_{\Sigma_2} \IK^p(w,h)\;  J(w,h) \; \alpha(h) \;d\sigma(w). 
\end{equation}
Then there is a subsequence converging uniformly on $\Sigma_2$ to a monotone surjection $h_\infty:\Sigma_2\to \Sigma_1$ which induces an isomorphism on fundamental groups.  The lift of $h_\infty$ to the universal covering space $\ID$,  $\tilde{h}_\infty:\overline{\ID}\to\overline{\ID}$ is a monotone surjection in $W^{1,2}(\ID)$ without infinite hairs, that is for all $\zeta\in \IS$,  $\tilde{h}_\infty^{-1}(\zeta)$ is a singleton in $\IS$. Further,  
\[\Phi= \IK^p(w,\tilde{h}) \, \tilde{h}_w(w) \, \overline{\tilde{h}_\wbar} \]
is holomorphic and automorphic.
\end{lemma}

The only thing not remarked on above is the statement about infinite hairs.  This follows since $\tilde{h}_\infty$ induces an isomorphism on fundamental groups.  An infinite hair lifts to a connected continuum in $\Sigma_2$.  This closure of this continuum,  call it $E$ (on which $h_\infty$ will be constant), can lie in no simply connected subset of $\Sigma_2$, since then it would lift into a single fundamental polyhedron and not go out to the boundary circle. Thus $E$ is topologically essential and since $h_\infty(E)$ is a single point,  this is impossible.  In fact note that this argument shows that any hair of $\tilde{h}_\infty$ has hyperbolic diameter bounded by that of $\Sigma_2$ (or what is the same thing,  a convex fundamental polyhedron for $\Gamma_2$.)

\medskip

We would now like to eliminate the possibility of hairs altogether.  In order to use our earlier ideas we need $\Phi\in L^1(\ID)$ which we know not to be the case.  Thus we localise. We also need to recall that $H=\tilde{h}_\infty$ is uniformly approximated by diffeomorphisms in $\overline{\ID}$.  Using this approximation we can see the following occurs.  Let $E$ be a hair, $H(E)=e\in\ID$. As $\Gamma_1(e)$ is a discrete set and $\gamma_1^i(e)\neq\gamma_1^j(e)$ if $i\neq j$,  the same must be true for $E$.  That is its translates under $\Gamma_2$ are disjoint.  Thus $E$ is compact and lifts into a simply connected open subset of some fundamental region $U$.  $U$ can be constructed from $H^{-1}(\ID(e,s))$ for $s$ sufficiently small.  Let
\[
V= H^{-1}(\ID(e,s)).
\]
Then $V$ is open,   as $H$ is continuous,  and as it is monotone one can see $U$ is simply connected.  Many of these facts about monotone mappings in two dimensions are collected in \cite{IO}.

Now let $X$ be a component of $int(\overline{V})\setminus V$.  These components are the cut points of $\partial V$.  As $H$ is continuous we have $H(X)\in\IS(e,r)$ and as $X$ is accessible from both sides we must have $H(X)$ a point.  That is $X$ is a hair in our previous terminology. Now $U=int(\overline{V})$ is also open and simply connected.   However we are going to have to use more theory to establish the following. (the hypothesis that $H$ is a homeomorphism outside the disk is simply to avoid possible pathologies at $\infty$ and can easily be improved.)

\begin{theorem}\label{thmcty} Let $H:\oD\to\oD$ be monotone.  Let $D=\ID(z_0,r)\subset\ID$,  $V=H^{-1}(D)$ and $
U =   int(\overline{V}).$ Then $U$ is simply connected and $\partial U$ has no cut points.  

Let $\varphi:\ID\to U$ be a Riemann map. Then the monotone map $H \circ \varphi:\ID\to\overline{D}$ admits a continuous extension to
\begin{equation}
H \circ \varphi : \overline{\ID} \to\overline{D} 
\end{equation}
which is monotone and  $H \circ \varphi(\IS)=\IS(r)$.
\end{theorem} 
\noindent{\bf Proof.} We have already discussed the first part here.  In what follows we retain our notation, so as to conform to the notation and definitions of \cite{Mil} and \cite{Rem}. 

\medskip

We need to consider the structure of the prime ends $\partial U$.  We use \cite{Rem} (note the interesting historical notes) and \cite[Chapter 17]{Mil} as modern accounts of the theory.  A quick review of the section ``Separation theorems and preliminaries'' of \cite{Rem} is recommended to those unfamiliar with the basics of prime end theory. The theory of prime ends was developed by Carath\'eodory 1913. Perhaps the main result is that the complement $K$ of a simply connected domain $U = \oC\setminus K$, $\#K > 1$, is locally connected if and only if every prime end has trivial impression. Equivalently, any Riemann map $\varphi : \ID \to U$ has a continuous extension to $\IS$. Let us briefly recall the definitions. A crosscut $C$ of a simply connected domain $U$ is a closed Jordan arc which intersects $\oC\setminus U$ exactly in its two endpoints. Every crosscut separates $U$ into precisely two components.

If $C$ is a crosscut of $U$ with a designated base point $b_0\in U\setminus C$, then $U\setminus C$ has exactly one component which does not contain $b_0$ denoted  $U_C$. A prime end of $U$ is represented by a sequence of pairwise disjoint crosscuts $\{C_n\}_{n\geq1}$ which satisfy the diameters $\diam(C_n)\to 0$ and $C_{n+1} \subset \overline{U_{C_n}}$. Two sequences $\{C_n\},  \{\tilde{C}_n\}$ represent the same prime end if  $\tilde{C}_j \subset U_{C_n}$ for all $n$ and all sufficiently large $j$, and vice versa.
The impression of a prime end $p$ is defined as
\begin{equation}
I(p)=\cap_{n} \overline{U_{C_n}} 
\end{equation}
and is a subset of $\partial U$. Note it is independent of the choice of representative sequence $\{C_n\}$ in the equivalence class. There is a natural way to define a topology on the set
\[ \hat{U} = U \cup \{p : p\; \mbox{is a prime end of $U$\}}. \]
A sequence $z_j\in U$ converges to a prime end $p$ if and only if $z_j\in U_{C_n}$ for all $n$ and all sufficiently large $j$. With this topology, $\hat{U}$  is homeomorphic to the closed unit disk $\overline{\ID}$. In fact the Riemann map extends continuously $\varphi : \overline{\ID} \to \hat{U}$, \cite[\S 3]{Rem}.

We now want to discuss what can be said if we have the additional piece of information that $H(\partial U)=\IS$ for a monotone map $H$.  First note that  $H|_{\partial U}:\partial U\to \IS(r)$ is monotone,  for if not there is a hair crossing $U$ (and so separating $U$). This easily leads to a contradiction as $H(V)=D$. Now we need the following lemma.

\medskip

\begin{lemma}\label{constant}
$H$ is constant on each prime end impression.
\end{lemma}

\medskip

We assume that Lemma \ref{constant} holds and complete the proof of Theorem \ref{thmcty} first. We have already noted that $\varphi : \overline{\ID} \to \hat{U}$ is a homeomorphism.  To see that $H\circ \varphi$ is monotone,  the only remaining issue is checking continuity at the impressions of prime ends.  But this follows from the definition of the topology on $\hat{U}$ and the fact that $H$ is uniformly continuous in a neighbourhood of $\partial U$ and constant on each impression.  We let $\varphi:\ID\to U$ be a Riemann map, so that $\varphi(0)\in E=H^{-1}(e)$ (recalling that $E$ is a noninfinite hair).  Then set
\begin{equation}
G = H \circ \varphi :\ID \to D.
\end{equation}
Note that $\IK(z,H)=\IK(\varphi(z),\tilde{h}_\infty)$ so that 
\begin{equation}
\Phi_G(w) = \IK^p(w,G) \, G_w(w) \, \overline{G_\wbar} = \Phi(\varphi)(\varphi')^2.
\end{equation}
It follows that $G:\ID\to D$ has holomorphic Ahlfors-Hopf differential in $L^1(\ID)$.  We know there is a solution without finite hairs,  and we are also in a position to assert uniqueness.  Thus $G$ is a homeomorphism,  and then a diffeomorphism.\hfill $\Box$

\medskip

\noindent{\bf Proof of Lemma \ref{constant}.} Let $p$ be a prime end and $I(p)$ its impression. We consider a representative sequence of crosscuts $C_n$. Here we may choose such a representative sequence that the base point $b_0$ has image $H(b_0)\notin \overline{H(U_{C_1})}\cup\partial D$. Indeed, for any given representative sequence of crosscuts $C_n$ of $I(p)$, there is an $N\geq1$ and a point $b_0'\in U\setminus\overline{U_{C_N}}$ such that $H(b_0')\notin \overline{H(U_{C_N})}\cup\partial D$, for otherwise
$\overline{D}\subset H(I(p))\cup\partial D$, which is impossible. Thus we may count from $N$ and replace the base point $b_0$ with $b_0'$ and this is the representative sequence we want. Note this also implies $H(b_0)\notin \overline{H(U_{C_n})}$ for all $n\geq1$ by definition, and then also $H(b_0)\notin \overline{H(C_n)}$ for all $n\geq1$.

We now consider the set $\overline{D}\setminus\overline{H(C_n)}$. As $H(b_0)\notin \overline{H(C_n)}$, there is a connected component $B_n$ of $\overline{D}\setminus\overline{H(C_n)}$ which contains $H(b_0)$. As $H$ is monotone, the set $H^{-1}(B_n)$ is a connected set in $\overline{U}$ that contains $b_0$ but no point in $C_n$. As $C_n$ separates $b_0$ and $U_{C_n}$ in $\overline{U}$, $H^{-1}(B_n)$ contains no point in $U_{C_n}$. We then have $H(U_{C_n})\subset\overline{D}\setminus B_n$. Now set $V_n:=\overline{D}\setminus\overline{B_n}$. Then $B_n$ and $V_n$ are two components which are separated by $C'_n:=\partial B_n\cap\overline{H(C_n)}$ in $\overline{D}$, and we also have $H(U_{C_n})\subset\overline{V_n}$. Since $H$ is continuous, $\diam(C_n)\to0$ implies $\diam(C_n')\to0$, and then as the image is the round disk $D$,
\[
\diam(H(U_{C_n}))\leq\diam(V_n)\to0.
\]
Thus
\[
H(I_p)\subset\lim_{n\to\infty}H(U_{C_n})\subset\lim_{n\to\infty}H(V_n),
\]
which is a single point.\hfill $\Box$

\medskip

\begin{center}
\scalebox{0.3}{\includegraphics*{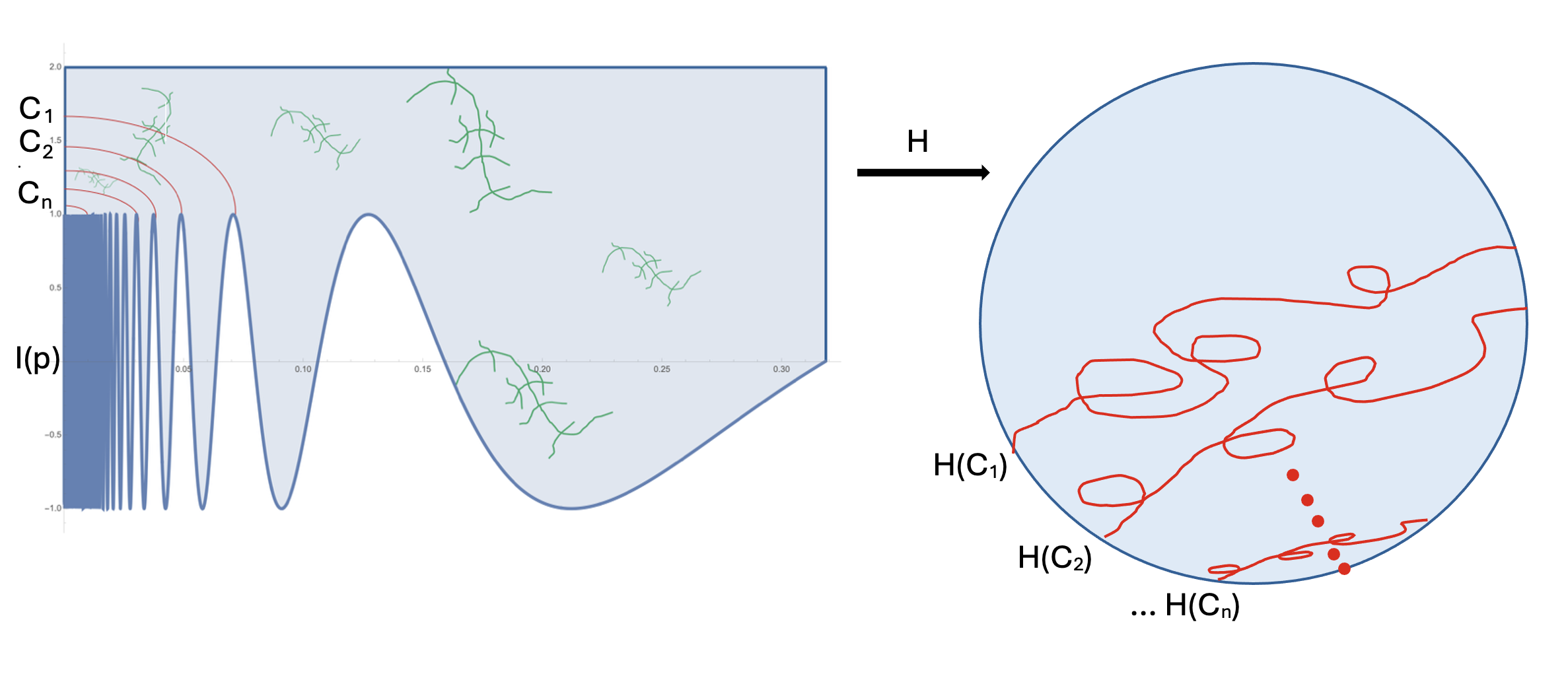}} \\
\end{center}
\vskip-15pt
\noindent {\bf Figure 5.} {\em An illustration of the preimage of a disk,  $H:\overline{U}\to \overline{D}$. The boundary contains a part $(x,-\sin \frac{\pi}{x})$, $x\in [-1,0]$ with prime end impression $[-i,i]$. The image of each $C_n$ might not be a path but must be a curve with finite length and separates $B_n$ from other parts of $\overline{D}\setminus\overline{H(C_n)}$ in $\overline{D}$. The map $H$ is uniformly continuous on $\overline{U}$, and as the diameter of $C_n\to 0$ we have $H(C_n)$ converging to a point.}

\medskip

This establishes the following.

\begin{theorem}  Given a pair of Riemann surfaces $\Sigma_1$ and $\Sigma_2$ and $p\geq 1$,  there is a unique minimiser to the problem
\[ \min_{f\in \, {\rm Diff}(\Sigma_1,\Sigma_2)} \int_{\Sigma_1} \IK^p(w,f) \; \alpha(w) \; d\sigma(w) \]
and this minimiser is a diffeomorphism with holomorphic Ahlfors-Hopf differential.
\end{theorem}

 \noindent{\bf Remark.} Actually our proof works for any weight $\alpha(w)$ as long as we can assert uniqueness in the local problem,  that is when the Ahlfors-Hopf differential is locally in $L^1(\ID)$. In particular $\alpha(w)\equiv 1.$

\end{document}